\documentclass[10pt]{article}
\usepackage{preamble}
\title{
\bf Error Estimation for Adaptive Mesh Refinement \\ in Droplet Simulations
}
\author{
Darsh Nathawani$^{1}$, Matthew G. Knepley$^{2}$ \\
$^{1}$Corresponding author: dnathawani@lsu.edu \\
$^{1}$Department of Civil and Environmental Engineering, Louisiana State University, \\ Baton Rouge, LA 70803, USA.\\
$^{2}$Department of Computer Science and Engineering, University at Buffalo, \\ Buffalo, NY 14260, USA.\\
}
\date{}

\begin{document}

\maketitle

{\bf Abstract:}
We present a one-dimensional shear-force-driven droplet formation model with a flux-based error estimator. The model is derived using asymptotic expansion and a front-tracking method to simulate the droplet interface. The model is then discretized using the Galerkin finite element method in the mixed form. However, the solution gradients exhibit large jumps across element boundaries and can grow rapidly due to the highly convective pinch-off process. This leads to an erroneous droplet interface and incorrect curvature. Therefore, the mesh must be sufficiently refined to capture the interface accurately. The mixed form of the governing equation naturally provides smooth interface gradients that can be used to compute the error estimate. The computed error estimate is then used to drive the adaptive mesh refinement algorithm. The efficacy of the error estimator is illustrated by comparing the droplet profiles obtained with adaptive refinement to those obtained with regular refinement. The adaptive mesh refinement approach reduces the computational cost significantly without compromising accuracy.

{\bf Keywords:} Droplet pinch-off, error estimation, adaptive mesh refinement, front-tracking method, mixed finite element method.


\section{Introduction}

The formation of droplets and bubbles is a mature field in physics. These surface-tension-driven phenomena have been deeply explored using interface-tracking and interface-capturing approaches. Historically, the accurate description of free surface flows emerged after Laplace~\cite{Laplace1805} and Young~\cite{Young1805} explained the role of mean curvature and the surface tension force as a source of instability in 1805. Eggers and Villermaux presented a comprehensive review of the liquid break-up process~\cite{EggersVillermaux2008}. Starting with the initial fluid column, the drop becomes heavier by continuously adding fluid. Under an external force like gravity, the instabilities start to grow on the interface. Meanwhile, the surface tension minimizes the surface energy by decreasing the radius of the fluid column. A spherical droplet starts to form at the end of the fluid column hanging from a very thin neck. Eventually, the radius becomes zero and a droplet separates from the initial fluid column. The point of zero radius is called the ``pinch-off" point. After the primary pinch-off, the neck recoils due to the unbalanced surface tension. The capillary waves due to this recoil perturb the interface before the tip can collapse back to the top of the fluid cone. This leads to a secondary breakup creating smaller droplets, which are called ``satellite drops''.

Inaugural work on the computational modeling of droplet formation was done by Eggers and Dupont~\cite{EggersDupont1994}. Their work involved finite difference formulation of the one-dimensional asymptotic model. Ambravaneswaran et al.~\cite{AmbravaneswaranWilkesBasaran2002} explored a finite element version of the same model. This model was only explored with the gravitational force as an external force. However, many industrial applications such as atomization and droplet entrainment in annular flow~\cite{Inoue_et_al_2021, BernaEscrivaMunozHerranz2015}, and microfluidic devices~\cite{Cramer2004, Teh2008, Dewandre2020microfluidic} involve droplet formation in a shear-force-driven environment. Using the same asymptotic approach for shear-induced droplets, we considered a one-dimensional model that estimates droplet sizes for an environment with co-flowing fluids~\cite{NathawaniKnepleyDropletShear}. For applications with a symmetrical structure of droplet formation, this one-dimensional model with the interface tracking approach is advantageous in terms of computational costs, while still maintaining the accuracy required for comparison with experiments. Nonetheless, this shear-induced droplet model deals with a wide range of external force magnitudes that induce a wide range of interface instability timescales.

The instability leading to the pinch-off point requires a finer mesh in the region of high curvature. As the interface evolves with the progression of numerical simulation, the mesh must be adaptively refined to control discretization errors, especially in the singularity regions. Therefore, a robust error estimation method is essential in designing an adaptive mesh refinement algorithm for our one-dimensional model. The field of \textit{a posteriori} error estimation has been under continuous and rigorous improvement since the 1980s~\cite{GratschBathe2005posteriori}. The preliminary work on error estimation was based on the finite element residuals and the flux jumps across the element boundaries~\cite{babuska1978posteriori, babuska1979analysis, babuska1987feedback, kelly1983posteriori}. Since then, many thorough review articles and books have been written on error estimation in finite element analysis; we refer the reader to~\cite{ainsworth1997posteriori, becker2001optimal, bangerth2003adaptive}. A thorough review of recent advancements in error estimation is given in~\cite{Chamoin_Legoll_Error_Estimation_Review}. In the present work, we show that the mixed finite element formulation of our one-dimensional model naturally provides an effective error estimator for adaptive mesh refinement.

In the following sections, we describe the one-dimensional droplet pinch-off model and explain the need for an error estimator. Then, we explore the flux-based error estimation approach and show that the mixed form naturally provides a smooth flux to compute error bounds. Finally, we illustrate the effectiveness of our error estimator using adaptive mesh refinement.

\section{Problem Description}

\begin{figure}[b!]
    \centering
    \includegraphics[width=0.3\textwidth]{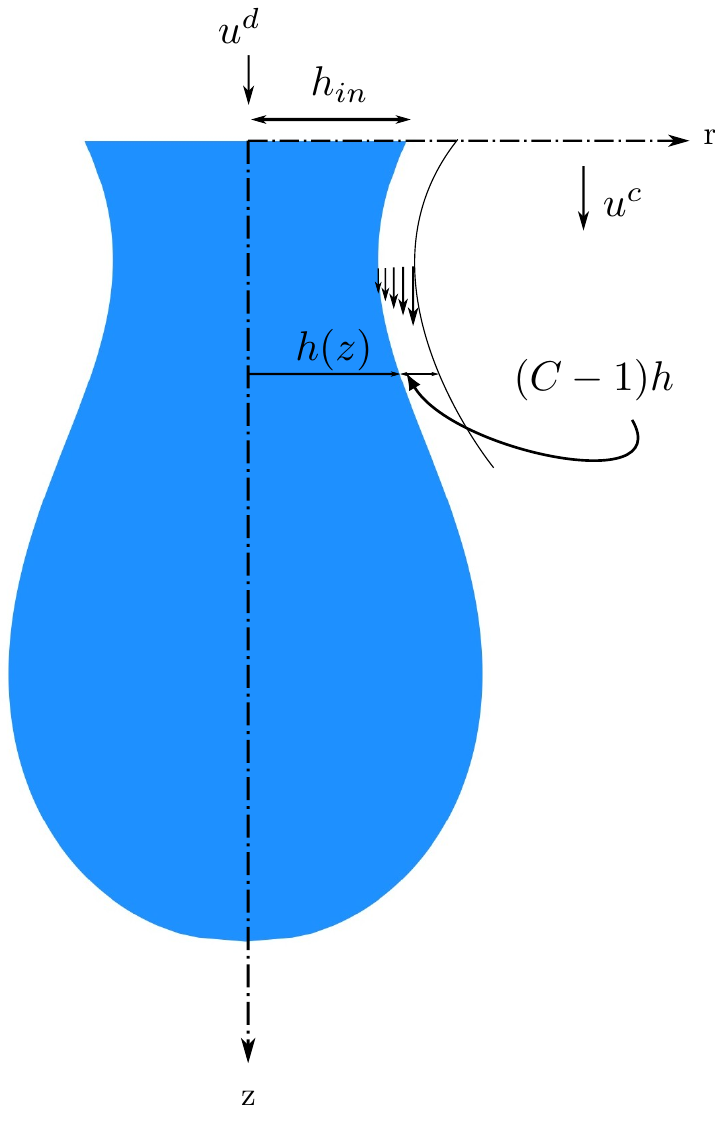}
    \caption{Schematic of a droplet formation in a co-flowing fluid.}
  \label{fig:coflow_schematic}
\end{figure}

We consider a fluid column slowly flowing with a constant velocity $u^d$ out of a nozzle with radius $h_{in}$. The \textit{continuous} (or outer) fluid is co-flowing with the \textit{dispersed} (or inner) fluid column but with a higher velocity $u^c$ in a capillary tube with radius $R$. Fig.~(\ref{fig:coflow_schematic}) shows a schematic of this co-flowing fluid scenario. The superscripts $d$ and $c$ represent dispersed phase and continuous phase fluids, respectively.

For the suspended fluid column without a background fluid, we start with the Navier-Stokes equations in cylindrical coordinates. The pinch-off process happens in finite time due to the surface tension forces trying to minimize the surface energy by contracting the interface in the radial direction. The radial contraction is faster than the axial elongation, which allows us to expand the solution variables asymptotically in $r$. Finally, the normal forces are balanced by the surface tension forces, which simplifies the equations to a one-dimensional model.

Similarly, we consider the droplet pinch-off of the fluid column in another co-flowing fluid under the gravitational and shear forces. We consider co-flowing fluids with symmetry preserved in the angular direction. This provides the basis for the one-dimensional approach. The continuous fluid affects the droplet morphology by applying shear forces on the interface between the two fluids. This shear-force effect can be approximated using an asymptotically derived force balance on the interface, yielding the one-dimensional governing equations~\cite{NathawaniKnepleyDropletShear}. The full derivation can be found in \cite{NathawaniPhDThesis}.

The momentum and the interface equations are given as follows.
\begin{align}
  \frac{\partial u^d}{\partial t} +  u^d \frac{\partial u^d}{\partial z} + \frac{\gamma}{\rho} \frac{\partial \mathcal{K}}{\partial z} - \frac{6 \nu^d}{h} \frac{\partial u^d}{\partial z} \frac{\partial h}{\partial z} \left( 1 + \frac{\mu^c}{\mu^d} \right)
  - 3 \nu^d \frac{\partial^2 u^d}{\partial z^2}\left( 1 + \frac{2}{3}\frac{\mu^c}{\mu^d} \right) &  \nonumber \\[2.5ex]
  + \frac{2}{\rho^d}\frac{dp^c}{dz} + \frac{1}{2 \rho^d \ln(C)}\frac{dp^c}{dz} - \left( 1 - \frac{\rho^c}{\rho^d} \right) g &= 0 \label{eq:co-flow_momentum} \\[3ex]
  \frac{\partial h}{\partial t} + u^d \frac{\partial h}{\partial z} + \frac{h}{2}\frac{\partial u^d}{\partial z} &= 0 \label{eq:interface}
\end{align}
Here, the parameters $\gamma$, $\rho$, $\nu$, $\mu$, and $p$ represent the surface tension coefficient, density, kinematic viscosity, dynamic viscosity, and pressure, respectively. The superscript represents dispersed or continuous phase fluid. $\mathcal{K}$ is the curvature defined by
\begin{align}
  \mathcal{K} = \left[ \frac{1}{h \left (1 + \frac{\partial h}{\partial z}^2\right )^{1/2}} - \frac{\frac{\partial^2 h }{\partial z^2}}{\left (1 + \frac{\partial h}{\partial z}^2\right )^{3/2}} \right] \label{eq:curvature}
\end{align}
The quantity $(C-1)h$ defines the thickness of the shear layer in the continuous phase flow, which defines how much force is experienced by the dispersed phase droplet. The parameter $C$ is a free parameter that can be defined using curve-fitted monotonic functions on the numerical and experimental data. In \cite{NathawaniKnepleyDropletShear}, it was shown that this single parameter model can account for behavior in a wide range of dimensionless parameter values.

The system of governing equations given by Eqs.~(\ref{eq:co-flow_momentum}-\ref{eq:interface}) is then discretized using a Galerkin finite element method in the mixed form~\cite{NathawaniKnepleyDropletShear}. We will use finite element spaces $\mathcal{U}$ for velocity, $\mathcal{H}$ for radius, and $\mathcal{S}$ for the slope variable, an auxiliary variable introduced in the mixed formulation to approximate the interface slope $\partial h / \partial z$. Now we can express the mixed finite element weak form as
\begin{align}
  &\int_{\Omega} q \left[ \frac{\partial u}{\partial t} +  u \frac{\partial u}{\partial z} -\frac{6\nu}{h}\frac{\partial h}{\partial z}\frac{\partial u}{\partial z}\left( 1 + \frac{\mu^c}{\mu^d} \right) + \frac{\gamma}{\rho} \left\{-\frac{s \frac{\partial s}{\partial z}}{h \left (1 + s^2\right )^{3/2}}  - \frac{s}{h^2 \left (1 + s^2\right )^{1/2}}\right\}  \right. \nonumber  \\[1.5ex]
  &\qquad \left.   + \frac{1}{2 \rho^d \ln(C)}\frac{dp^c}{dz} - \left( 1 - \frac{\rho^c}{\rho^d} \right)g \right] d\Omega  + \int_{\Omega} \nabla q \left[3 \nu \left( 1 + \frac{2}{3}\frac{\mu^c}{\mu^d} \right) \frac{\partial u}{\partial z}\right.  \nonumber\\[2ex]
  &\qquad \left. + \frac{\gamma}{\rho} \frac{\frac{\partial s}{\partial z}}{\left (1 + s^2\right )^{3/2}}  \right] d \Omega   - \int_{\Gamma} q \left[3 \nu \left( 1 + \frac{2}{3}\frac{\mu^c}{\mu^d} \right) \frac{\partial u}{\partial z} + \frac{\gamma}{\rho} \frac{\frac{\partial s}{\partial z}}{\left (1 + s^2\right )^{3/2}}  \right] d \Gamma = 0   \label{eq:FE_momentum}  \\[3ex]
  &\int_{\Omega} v \left [ \frac{\partial h}{\partial t} + u \frac{\partial h}{\partial z} + \frac{h}{2} \frac{\partial u}{\partial z} \right ]d \Omega = 0 \label{eq:FE_interface} \\[3ex]
 &\int_{\Omega} w \left [ s - \frac{\partial h}{\partial z} \right ] d \Omega = 0  \label{eq:FE_s}
\end{align}
where $u, q \in \mathcal{U}$, $h, v \in \mathcal{H}$, and $s, w \in \mathcal{S}$. Eqs.~(\ref{eq:FE_momentum}-\ref{eq:FE_s}) are solved using the continuous Galerkin method with $C^0(\Omega)$ elements. The solution algorithm involves a moving mesh and calculates the droplet length in a self-consistent way. The algorithm is explained in \cite{NathawaniKnepleyDropletGravity}, and has been implemented using the PETSc libraries~\cite{petsc-user-ref,petsc-web-page}.

The flow becomes highly convective as the droplet interface evolves up to the point at which the neck begins to form. The radial contraction of the interface accelerates faster than the axial elongation, rapidly approaching the singularity. This behavior must be captured accurately to obtain precise equilibrium droplet profiles. In addition, the interface is advected with the flow, so that errors in the interface location at one step strongly affect later steps in the flow.

\section{Error estimation}

When approaching the singularity, the high-curvature regions quickly develop high stresses. An error in the solution can cause the interface to advect in an incorrect direction, yielding an incorrect droplet profile. A coarser mesh must be refined to improve the solution in those regions and obtain an accurate force balance, leading to more accurate droplet profiles. This refinement can be done adaptively by estimating the true error in the solution using an \textit{a posteriori} error estimation approach. The foremost goal is to compute the discretization error in the droplet interface evolution.

Since the solution is approximated with $C^0(\Omega)$ continuity, the solution gradients are discontinuous across the element boundaries. This discontinuity significantly impacts the curvature computation and can produce erroneous equilibrium droplet profiles. We use a flux-recovery-based approach to drive the adaptive mesh refinement. The flux-recovery approach aims to post-process a smooth gradient from the finite element solution and construct the error estimator from the difference between the smoothed and non-smoothed gradients.

In the droplet pinch-off model, the solution consists of the velocity ($u$) and the droplet radius ($h$). The mixed formulation additionally includes an approximation of the gradient $\partial h/\partial z$ as part of the solution, which is denoted by $s$. In other words, the mixed variable $s$ already provides a smooth gradient of $h$. If the true gradient is $\partial \bar{h} / \partial z$, the quantity $s$ is assumed to retain better accuracy than $\partial h/\partial z$ since the slope (or flux) $s$ is part of the equilibrium equations and coupled with the curvature derivatives in the momentum equations. The true error norm in the gradient can be written as
\begin{align}\label{eq:compErrorBound}
  \left|\left| e \right|\right| &= \left|\left| \frac{\partial \bar{h}}{\partial z} - \frac{\partial h}{\partial z} \right|\right| \approx \left|\left| s - \frac{\partial h}{\partial z} \right|\right|
\end{align}
The quality of the estimate clearly depends on how good the approximation $s$ is of the true slope $\partial \bar{h} / \partial z$. Assuming there exists $c \in (0,1)$ such that
\begin{align}
  \left|\left| \frac{\partial \bar{h}}{\partial z} - s \right|\right| \leq c  \left|\left| \frac{\partial \bar{h}}{\partial z} - \frac{\partial h}{\partial z} \right|\right|
\end{align}
Using the triangle inequality, the bounds on the true error in terms of the error estimator follow directly:
\begin{align}
  \frac{\eta}{1 + c} \leq \left|\left| e \right|\right| \leq \frac{\eta}{1 - c} \label{eq:error-bounds}
\end{align}
where $\eta =  \left|\left| s - \partial h / \partial z \right|\right| $ is the error estimate of the true error. It is clear that the smaller the $c$, the better the error estimate.

The traditional \textit{a posteriori} error estimation methods, such as residual-based and flux-based methods, use the residual in each element and the flux jumps across the elemental interfaces. However, these error bounds include an unknown multiplicative constant. Our approach can also be viewed from the residual-based error estimation perspective. A residual-based a posteriori error estimation approach directly utilizes the finite element solution. Defining the error in the energy norm, one can also derive upper and lower error bounds for the error in the slope $s$ using the element residuals of Eq.~(\ref{eq:FE_s}).

Here, the error bound given in Eq.~(\ref{eq:compErrorBound}) are used to refine the mesh adaptively. A simple $h$-refinement can be done cyclically:
\begin{center}
  solve $\rightarrow$ estimate error $\rightarrow$ mark elements $\rightarrow$ refine
\end{center}

\noindent We use two strategies to refine the mesh and compare their effectiveness.

\noindent(1) The regular refinement, where elements are doubled at each refinement cycle.

\noindent(2) The elements are marked using D\"{o}rfler's method~\cite{dorfler1996}, where the elements are marked such that the total error in the marked elements is $\lambda$ times the total error in the entire domain.

\section{Results}

The droplet simulations presented in this section are for an 85\% glycerol solution with an inlet velocity of 5 mm/s and a co-flowing air stream at 1 m/s. The nozzle inlet radius is 2.5 mm and the outer capillary tube radius is 25 mm. The material properties are taken from~\cite{NathawaniKnepleyDropletGravity}.

\begin{figure}[h]
    \centering
    \begin{subfigure}[t]{0.2\textwidth}
        \centering
        \includegraphics[width=0.8\linewidth, valign=t]{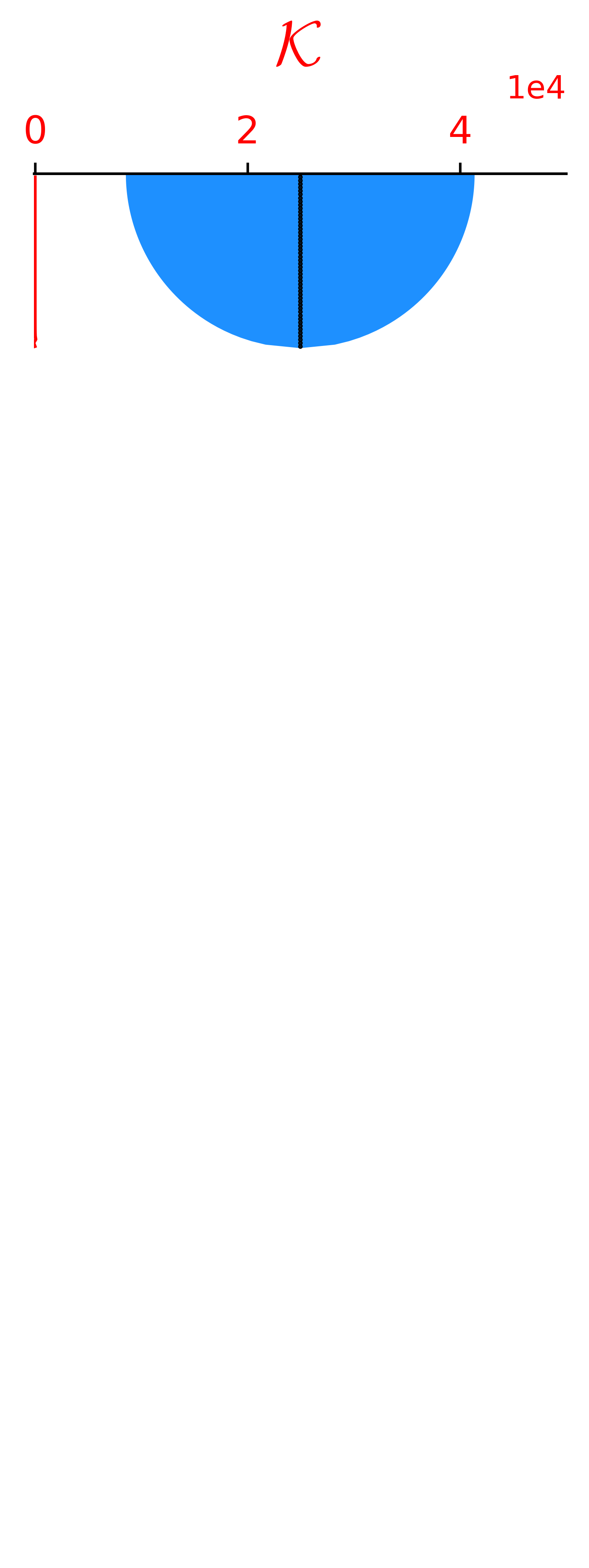}
    \end{subfigure}\hfill%
    \begin{subfigure}[t]{0.2\textwidth}
        \centering
        \includegraphics[width=0.8\linewidth, valign=t]{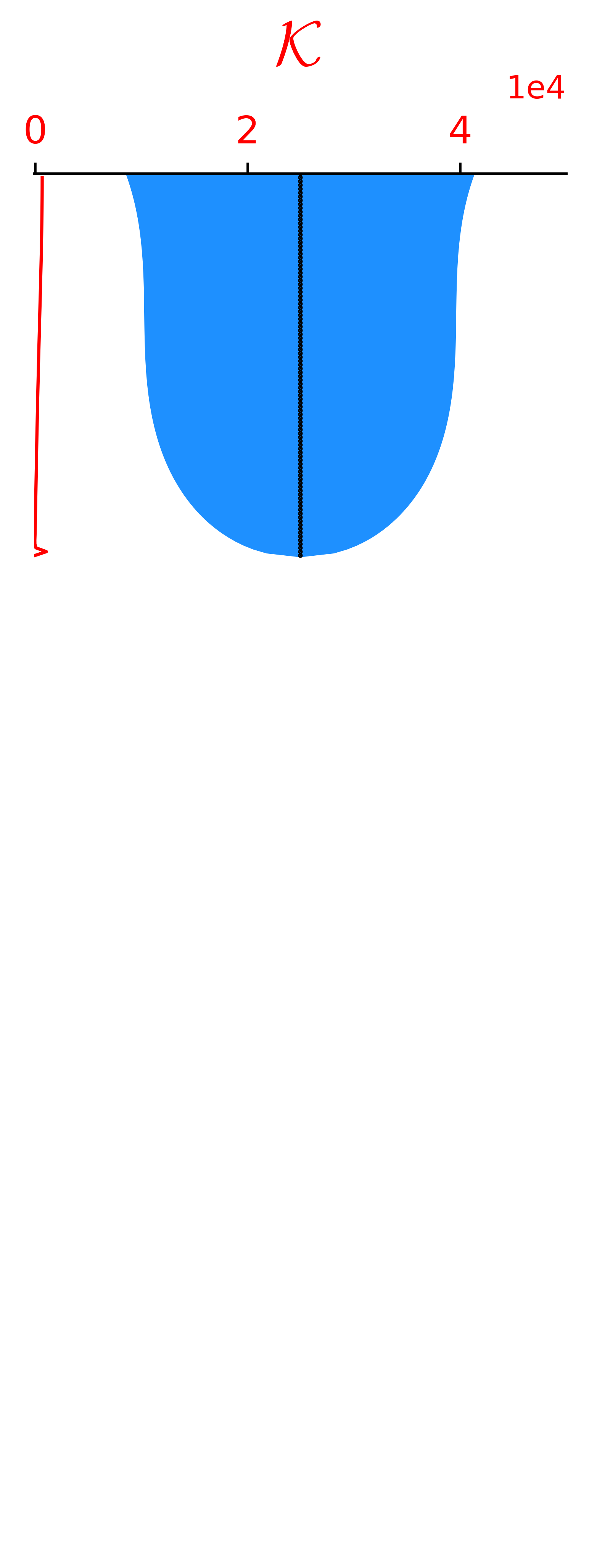}
    \end{subfigure}\hfill%
    \begin{subfigure}[t]{0.2\textwidth}
        \centering
        \includegraphics[width=0.8\linewidth, valign=t]{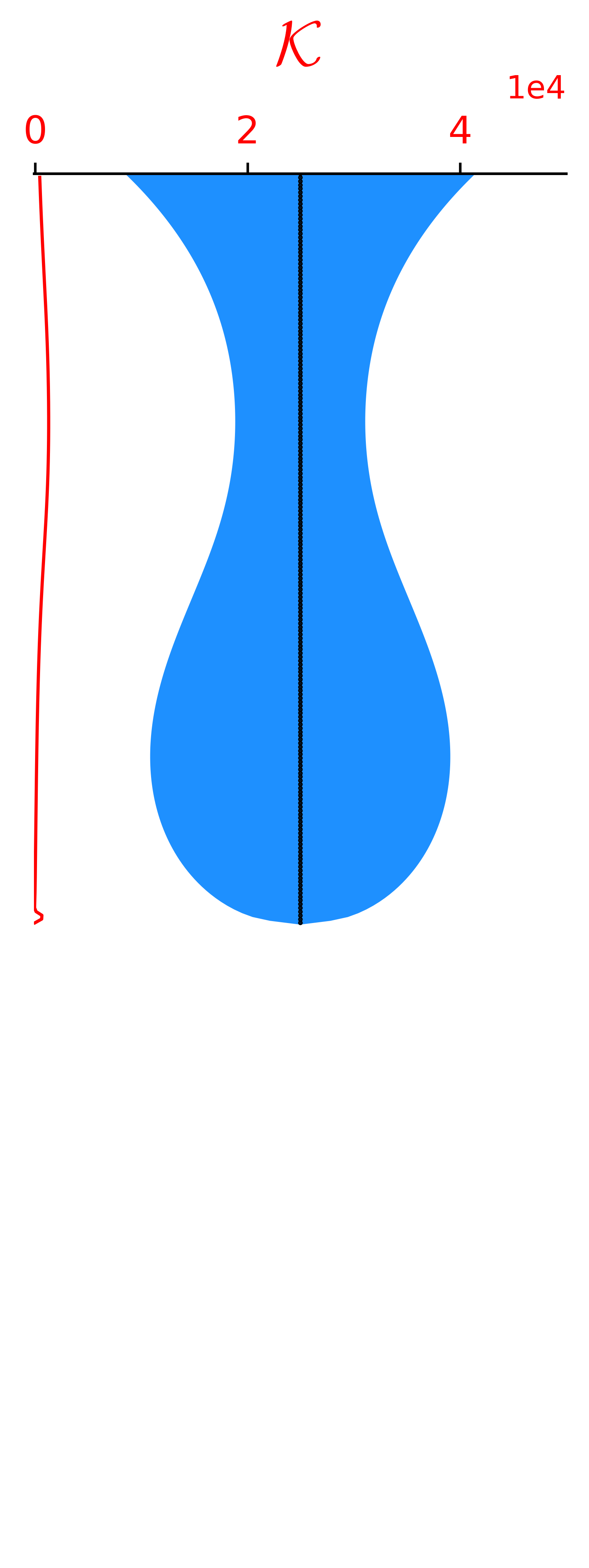}
    \end{subfigure}\hfill%
    \begin{subfigure}[t]{0.2\textwidth}
        \centering
        \includegraphics[width=0.8\linewidth, valign=t]{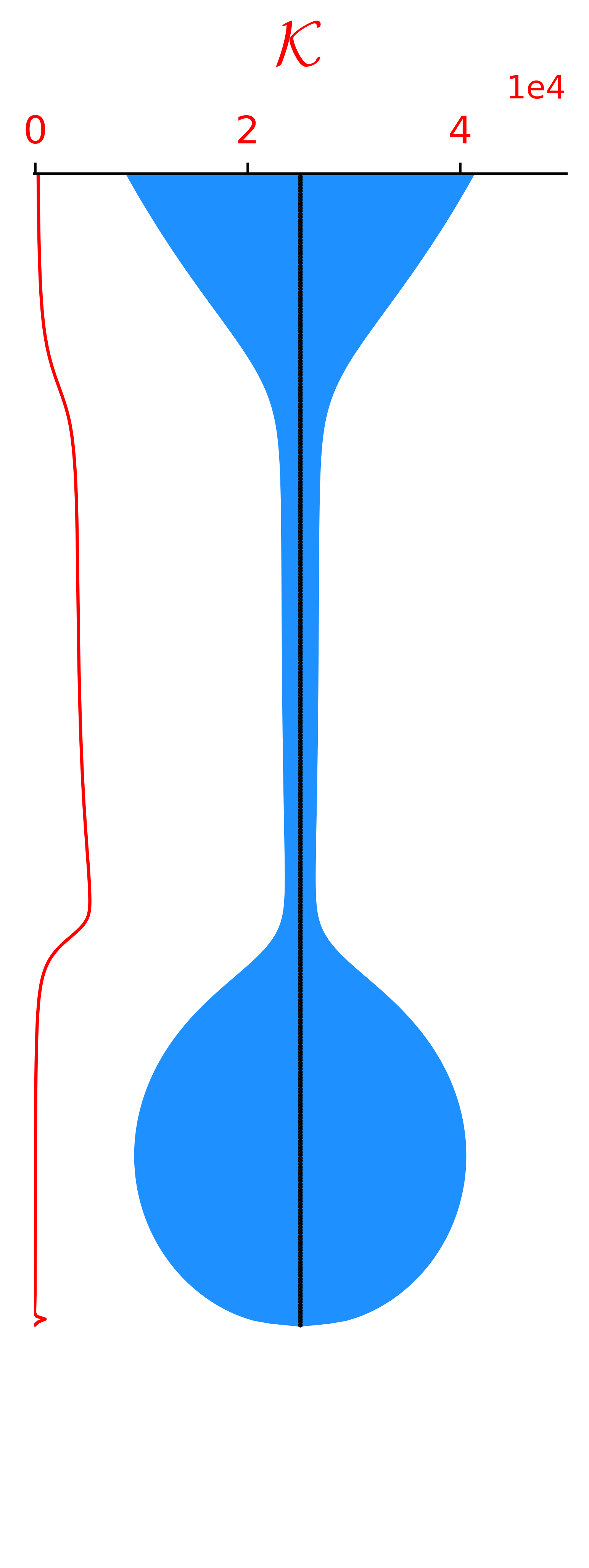}
    \end{subfigure}\hfill%
    \begin{subfigure}[t]{0.2\textwidth}
        \centering
        \includegraphics[width=0.8\linewidth, valign=t]{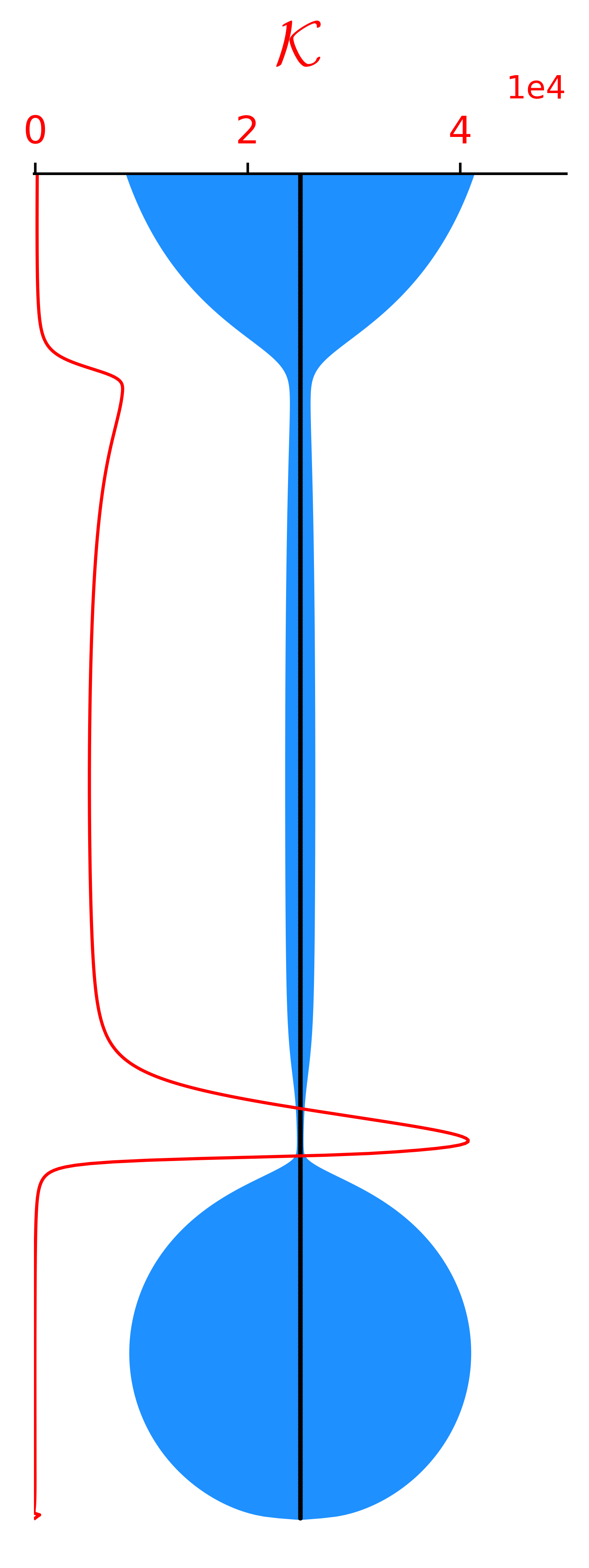}
    \end{subfigure}
    \caption{Change in the mean curvature with the droplet evolution.}
    \label{fig:curvature-change}
\end{figure}

The error estimation approach is motivated entirely by the need to capture regions with sharp changes in the mean curvature. Therefore, it is important to understand how the mean curvature evolves in order to contextualize the error estimates. Fig.~(\ref{fig:curvature-change}) shows a sequence of images illustrating the evolution of mean curvature as the droplet evolves up to the primary pinch-off. The mean curvature is computed using the finite element solution and plotted along the droplet length. The initial mesh is uniform with 50 elements. As the droplet length increases, the mesh is regularly refined by doubling the number of elements as the length increases by a factor of 2 to maintain the accuracy of the solution. The droplet evolves roughly up to 8 times the initial length before the primary pinch-off. Therefore, the mesh is refined up to 800 elements, which 4 levels of refinements. The refinement cycle is defined by the non-dimensionalized droplet length $L/h_{in}$, where $L$ is the droplet length and $h_{in}$ is the inlet radius. We refine the mesh when $L/h_{in} > N$, where $N = \{2, 4, 6\}$ for the first, second, and third refinement cycle. The last refinement cycle is triggered when $h_{min} < \frac{h_{in}}{50}$, where $h_{min}$ is the minimum droplet radius. This ensures that the mesh is sufficiently refined to capture the pinch-off process. The pinch-off is concluded when $h_{min} < \frac{h_{in}}{100}$.

Initially, the mean curvature is largest at the droplet tip. The mean curvature starts increasing in the neck region of the fluid column as the droplet evolves. This is naturally due to the radial shrinkage of the column. Finally, the mean curvature is highest at the pinch-off location when the primary droplet separates from the fluid column. The regions with sharp curvature changes are the top of the neck, the region near the pinch-off, and the droplet tip. Therefore, a higher error is expected in these regions. One notable observation is that the curvature has a kink at the tip of the droplet. This is due to the choice of initial value of $s = -10$ at the tip and kept bounded with a large negative value. However, the solution is invariant to the choice of this initial value as long as it is sufficiently negative. Theoretically, the slope at the tip should be negative infinity. The choice of a large negative value is a practical way to approximate this behavior and to keep the simulation stable.

\begin{figure}[h]
    \centering
    \begin{subfigure}[t]{0.48\textwidth}
        \centering
        \includegraphics[width=0.35\linewidth]{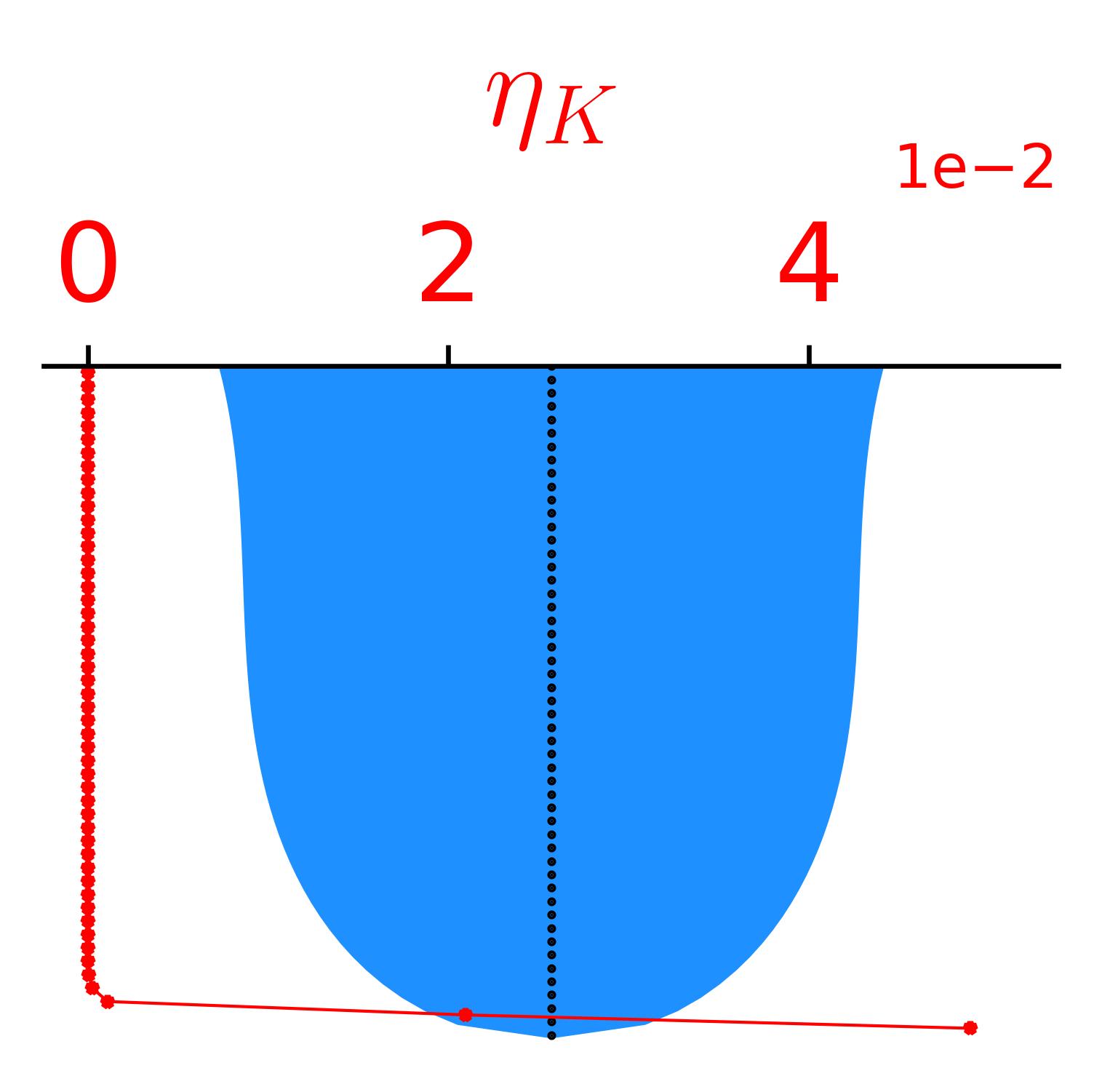}
        \caption{Initial droplet phase.}
        \label{fig:initial-droplet-error}
    \end{subfigure}
    \begin{subfigure}[t]{0.48\textwidth}
        \centering
        \begin{tikzpicture}
            \node[inner sep=0pt] (img) at (0,0) {\includegraphics[width=0.35\linewidth]{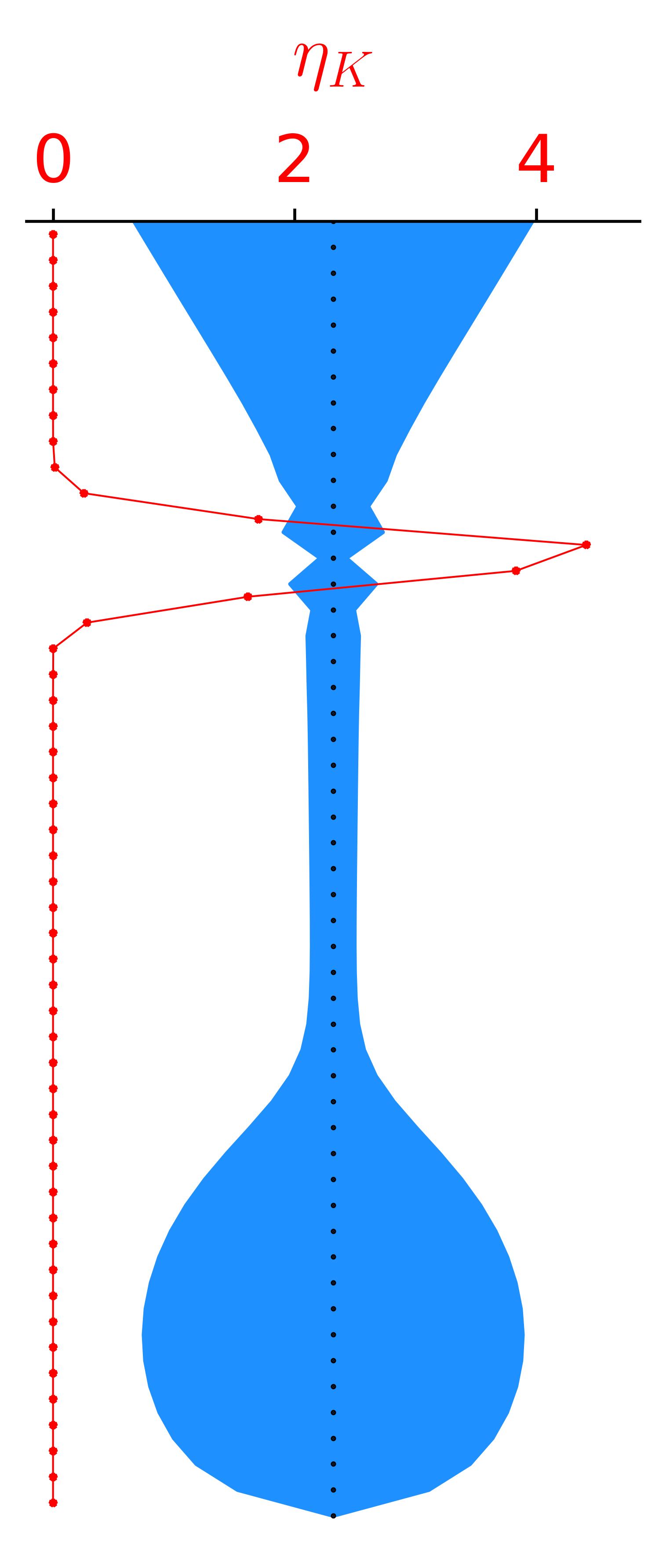}};
            \node[draw, green, minimum height=2.1cm, minimum width=3.0cm] (box) at (-0.,-2.0) {};
            \node[inner sep=0pt] (img-cut) at (4,0) {\includegraphics[width=0.35\linewidth]{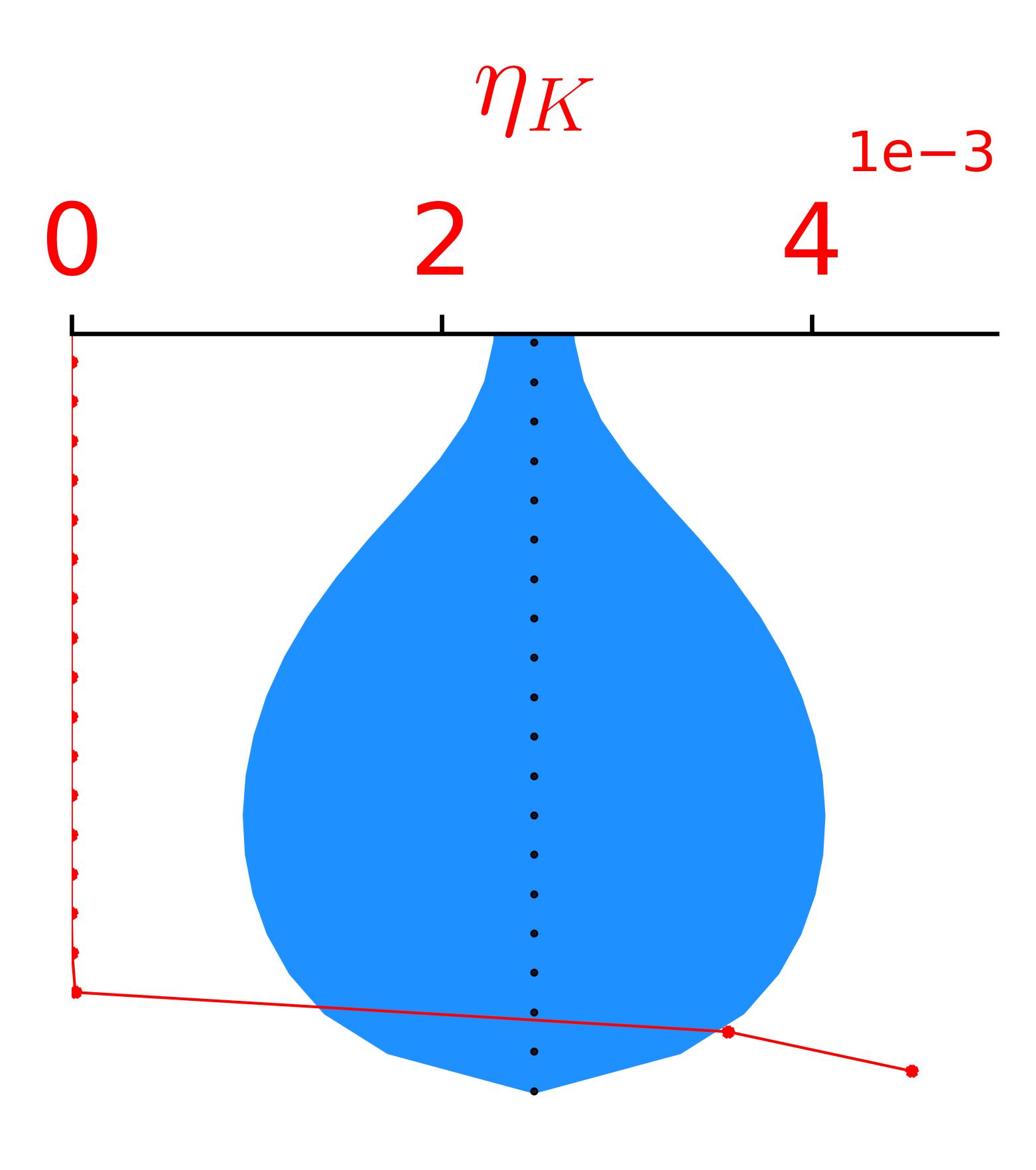}};
            \draw[-stealth, green] (box.east) -- (img-cut);
        \end{tikzpicture}
        \caption{Evolved droplet phase with a magnified view at the tip.}
        \label{fig:evolved-droplet-error-zoom}
    \end{subfigure}
    \caption{Error plots showing the error magnitude at the initial phase and the evolved phase.}
    \label{fig:error-initial-evolved}
\end{figure}

We first report the element-wise errors without performing adaptive mesh refinement. The droplet length and radius are non-dimensionalized using the inlet radius $h_{in}$ in all plots. The element-wise estimated error, denoted by $\eta_K$, is shown in Fig.~(\ref{fig:error-initial-evolved}). The error is plotted for the initial droplet profile and an evolved droplet profile as shown in Fig.~(\ref{fig:initial-droplet-error}) and Fig.~(\ref{fig:evolved-droplet-error-zoom}), respectively. The initial droplet profile has the highest curvature gradient near the tip. Thus, the highest error is detected at the tip, as shown in the figure.

On the other hand, the evolved droplet has a large error at the top of the neck region. This error is associated with an inaccurate interface representation in the neck region. This inaccuracy leads to a divergence of residuals and simulation failure before reaching the primary pinch-off. The magnified image shows the error magnitude at the tip, which is larger for the last few elements. This suggests that mesh refinement is required to accurately capture the smooth spherical shape at the droplet tip.

\begin{figure}[h!]
\centering
\begin{subfigure}[t]{0.45\textwidth}
    \centering
    \includegraphics[width=0.5\linewidth]{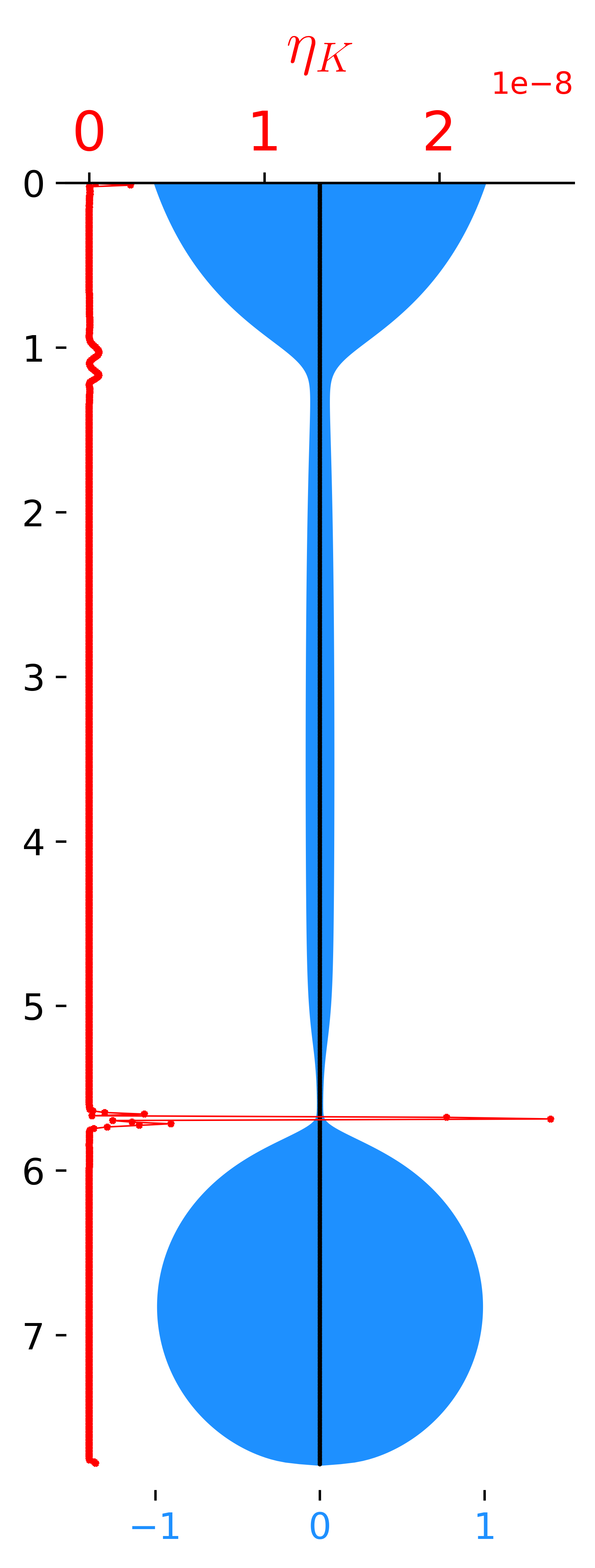}
    \caption{Regular refinement.}
  \label{fig:fine-error-plot-regular}
\end{subfigure}
\hfill
\begin{subfigure}[t]{0.45\textwidth}
    \centering
    \includegraphics[width=0.5\linewidth]{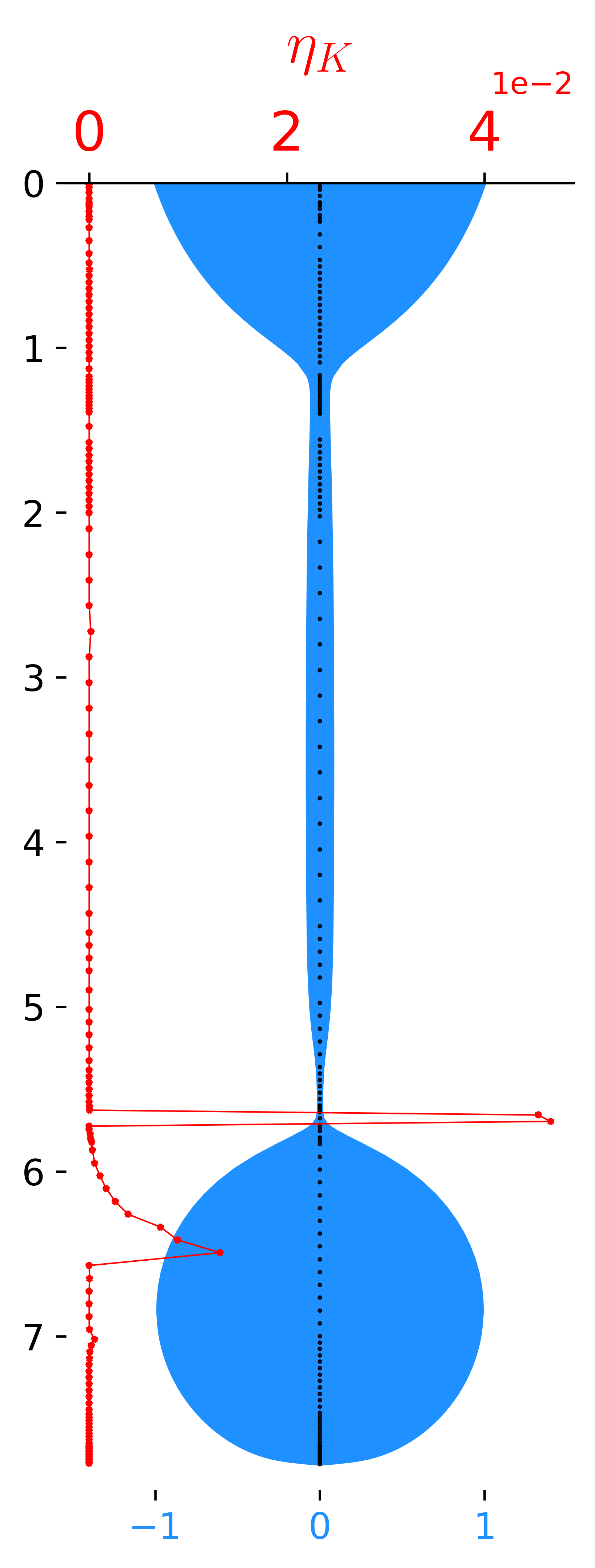}
    \caption{Adaptive refinement using D\"{o}rfler strategy.}
  \label{fig:fine-error-plot-dorfler}
\end{subfigure}
\caption{Comparison of error plots superposed on the 85\% glycerol droplet.}
\label{fig:error-plots}
\end{figure}

Based on the reported errors, the adaptive mesh refinement strategy described in the previous section is applied. As the simulation progresses, we track the error evolution along with the droplet interface. This provides insight into how the error evolves as the curvature changes. Using the  D\"{o}rfler's strategy, adaptive mesh refinement is performed when the droplet length is greater than $N$ times the inlet radius, where $N = \{2,4,6,...\}$. Similar to the regular refinement strategy, the final refinement is performed when $h_{min} < \frac{h_{in}}{50}$. For the D\"{o}rfler's approach, the elements with a cumulative error equal to 30\% of the total error are marked and refined. The choice of these parameters is based on observations from the error evolution on a unrefined grid.

Fig.~(\ref{fig:error-plots}) directly compares regular refinement (left) and adaptive refinement using D\"{o}rfler marking (right) at the same evolution stage. The regular strategy, shown in Fig.~(\ref{fig:fine-error-plot-regular}), reduces the error by uniformly adding elements over the full domain. As a result, the error in the neck region is significantly reduced. The largest error is at the pinch-off location, where the curvature changes are highest. Additionally, the error at the top of the neck and at the tip is also apparent. Nonetheless, the error is reduced to the order of $10^{-8}$, which is sufficient to capture the droplet profile accurately. We use this as a reference to compare the D\"{o}rfler marking approach, which refines only the elements with the largest error.

Fig.~(\ref{fig:fine-error-plot-dorfler}) shows the error distribution for the D\"{o}rfler marking approach. Similar to the regular refinement, the error is largest at the pinch-off location. The error is reduced to the order of $10^{-2}$, which is smaller than the unrefined case but larger than the regular refinement case. However, the droplet profile is significantly improved compared to the unrefined case, and the instabilities are removed. Moreover, the pinch-off droplet profile is visually indistinguishable from the regular refinement case, indicating that the error level is sufficient to capture the droplet profile accurately.

Comparing Fig.~(\ref{fig:fine-error-plot-regular}) and Fig.~(\ref{fig:fine-error-plot-dorfler}), both approaches recover a smooth interface and remove the visible instability near the neck. However, a quantitative comparison is also necessary to fully assess the performance of the D\"{o}rfler marking approach. Table~(\ref{tab:regular-vs-amr-qoi}) compares the key quantities of interest and simulation time for a representative regular-refinement run and an AMR run. The quantities are pinch-off location, surface area, volume, and pinch-off time. The AMR case reproduces all quantities with sub $1\%$ inaccuracy compared to the regular-refinement reference, except the pinch-off time which has $1.3174\%$ difference. This is achieved while reducing wall-clock time from $638$ s to $153$ s, which is $76\%$ reduction in computational cost. In other words, this AMR provides a speedup of $4.17\times$, while maintaining close agreement in all quantities of interest. This comparison is evident of the effectiveness of the error estimator in driving the adaptive mesh refinement to capture the droplet profile accurately while significantly reducing computational cost.

The key point is that AMR achieves nearly the same quantities of interest while using far fewer elements. The regular-refinement case reaches 800 elements, whereas the AMR case uses at most 146 elements. This large reduction in mesh size is the main reason for the substantial savings in computational cost.

\begin{table}[h!]
\centering
\small
\caption{Comparison of regular refinement and AMR for pinch-off metrics and computational cost. Relative differences are computed with respect to the regular-refinement result.}
\label{tab:regular-vs-amr-qoi}
{\setlength{\tabcolsep}{9pt}
\renewcommand{\arraystretch}{1.2}
\begin{tabular}{lccc}
\hline
Quantity & Regular refinement & AMR & Relative difference \\
\hline
Non-dimensional pinch-off location & 0.7200 & 0.7186 & 0.1944\% \\
Pinch-off droplet surface area & $7.5548 \times 10^{-5}$ & $7.5964 \times 10^{-5}$ & 0.5506\% \\
Pinch-off volume & $6.1558 \times 10^{-8}$ & $6.2043 \times 10^{-8}$ & 0.7878\% \\
Pinch-off time (s) & 0.5769 & 0.5845 & 1.3174\% \\
Wall-clock time (s) & 638 & 153 & 76.02\% reduction \\
Maximum number of elements & 800 & 146 & 81.75\% reduction \\
\hline
\end{tabular}
}
\end{table}

The one-dimensional model has been previously used for uncertainty quantification (UQ) analysis of droplet atomization in hybrid rocket combustion~\cite{Georgalis2024UQ}. The error estimator can be used to drive the adaptive mesh refinement in the UQ simulations, allowing for more accurate and efficient exploration of the parameter space. This is particularly important in UQ analyses, where a large number of simulations are required to capture the variability in the system. By using the error estimator to refine the mesh adaptively, we can ensure that each simulation is performed with sufficient accuracy while minimizing computational costs.

\section{Conclusion}

We presented a simple yet effective strategy to compute the error estimate for droplet simulations with a front-tracking method. The one-dimensional model used for the simulations is accurate for shear-induced droplets in a co-flowing environment. The mixed form equations naturally incorporate a smooth flux as part of the solution, enabling flux-based error estimation. The comparison of error evolution for 85\% glycerol droplets in co-flowing air is presented for unrefined, regularly refined, and adaptively refined grids. This comparison shows that the adaptively refined grid provides improved stability and accuracy in capturing the droplet interface. The AMR approach keeps the quantities of interest within $1\%$ of the regular-refinement case while reducing the computational cost by $\approx 76\%$ by using $\approx 82\%$ fewer elements.

\printbibliography

@article{Chamoin_Legoll_Error_Estimation_Review,
  title={An introductory review on a posteriori error estimation in finite element computations},
  author={Chamoin, Ludovic and Legoll, Fr{\'e}d{\'e}ric},
  journal={SIAM Review},
  volume={65},
  number={4},
  pages={963--1028},
  year={2023},
  publisher={SIAM}
}

@article{NathawaniKnepleyDropletGravity,
author = {Nathawani, Darsh K. and Knepley, Matthew G. },
title = {Droplet formation simulation using mixed finite elements},
journal = {Physics of Fluids},
volume = {34},
number = {6},
pages = {064105},
year = {2022},
doi = {10.1063/5.0089752}
}

@article{NathawaniKnepleyDropletShear,
  title={A one-dimensional mathematical model for shear-induced droplet formation in co-flowing fluids},
  author={Nathawani, Darsh and Knepley, Matthew},
  journal={Theoretical and Computational Fluid Dynamics},
  pages={1--17},
  year={2024},
  publisher={Springer}
}

@phdthesis{NathawaniPhDThesis,
  title={Droplet Formation: One-Dimensional Mathematical Model and Computations},
  author={Nathawani, Darsh K},
  year={2023},
  school={State University of New York at Buffalo}
}

@article{EggersDupont1994,
  title={Drop formation in a one-dimensional approximation of the Navier--Stokes equation},
  author={Eggers, Jens and Dupont, Todd F},
  journal={Journal of fluid mechanics},
  volume={262},
  pages={205--221},
  year={1994},
  publisher={Cambridge University Press}
}

@article{AmbravaneswaranWilkesBasaran2002,
  title={Drop formation from a capillary tube: Comparison of one-dimensional and two-dimensional analyses and occurrence of satellite drops},
  author={Ambravaneswaran, Bala and Wilkes, Edward D and Basaran, Osman A},
  journal={Physics of Fluids},
  volume={14},
  number={8},
  pages={2606--2621},
  year={2002},
  publisher={American Institute of Physics}
}

@article{GratschBathe2005posteriori,
  title={A posteriori error estimation techniques in practical finite element analysis},
  author={Gr{\"a}tsch, Thomas and Bathe, Klaus-J{\"u}rgen},
  journal={Computers \& structures},
  volume={83},
  number={4-5},
  pages={235--265},
  year={2005},
  publisher={Elsevier}
}

@inproceedings{Georgalis2024UQ,
  title={Uncertainty Quantification of Shear-induced Paraffin Droplet Pinch-off in Hybrid Rocket Motors},
  author={Georgalis, Georgios and Nathawani, Darsh and Knepley, Matthew and Patra, Abani},
  booktitle={AIAA SCITECH 2024 Forum},
  pages={1021},
  year={2024}
}

@article{Laplace1805,
  title={Trait{\'e} de m{\'e}canique c{\'e}leste, vol. 4},
  author={Laplace, Pierre Simon},
  journal={Supplements au Livre X},
  year={1805}
}

@article{Young1805,
  title={III. An essay on the cohesion of fluids},
  author={Young, Thomas},
  journal={Philosophical transactions of the royal society of London},
  number={95},
  pages={65--87},
  year={1805},
  publisher={The Royal Society London}
}

@article{EggersVillermaux2008,
  title={Physics of liquid jets},
  author={Eggers, Jens and Villermaux, Emmanuel},
  journal={Reports on progress in physics},
  volume={71},
  number={3},
  pages={036601},
  year={2008},
  publisher={IOP Publishing}
}

@article{Cramer2004,
  title={Drop formation in a co-flowing ambient fluid},
  author={Cramer, Carsten and Fischer, Peter and Windhab, Erich J},
  journal={Chemical Engineering Science},
  volume={59},
  number={15},
  pages={3045--3058},
  year={2004},
  publisher={Elsevier}
}

@article{Inoue_et_al_2021,
  title={On the droplet entrainment from gas-sheared liquid film},
  author={Inoue, Chihiro and Maeda, Ikkan},
  journal={Physics of Fluids},
  volume={33},
  number={1},
  pages={011705},
  year={2021},
  publisher={AIP Publishing LLC}
}

@article{BernaEscrivaMunozHerranz2015,
  title={Review of droplet entrainment in annular flow: Characterization of the entrained droplets},
  author={Berna, C and Escriv{\'a}, A and Mu{\~n}oz-Cobo, JL and Herranz, LE},
  journal={Progress in Nuclear Energy},
  volume={79},
  pages={64--86},
  year={2015},
  publisher={Elsevier}
}

@article{Dewandre2020microfluidic,
  title={Microfluidic droplet generation based on non-embedded co-flow-focusing using 3D printed nozzle},
  author={Dewandre, Adrien and Rivero-Rodriguez, Javier and Vitry, Youen and Sobac, Benjamin and Scheid, Benoit},
  journal={Scientific reports},
  volume={10},
  number={1},
  pages={1--17},
  year={2020},
  publisher={Nature Publishing Group}
}

@article{Teh2008,
  title={Droplet microfluidics},
  author={Teh, Shia-Yen and Lin, Robert and Hung, Lung-Hsin and Lee, Abraham P},
  journal={Lab on a Chip},
  volume={8},
  number={2},
  pages={198--220},
  year={2008},
  publisher={Royal Society of Chemistry}
}

@article{babuska1978posteriori,
  title={A-posteriori error estimates for the finite element method},
  author={Babu{\v{s}}ka, Ivo and Rheinboldt, Werner C},
  journal={International journal for numerical methods in engineering},
  volume={12},
  number={10},
  pages={1597--1615},
  year={1978},
  publisher={Wiley Online Library}
}

@article{babuska1987feedback,
  title={A feedback finite element method with a posteriori error estimation: Part I. The finite element method and some basic properties of the a posteriori error estimator},
  author={Babu{\v{s}}ka, Ivo and Miller, A},
  journal={Computer Methods in Applied Mechanics and Engineering},
  volume={61},
  number={1},
  pages={1--40},
  year={1987},
  publisher={Elsevier}
}

@article{babuska1979analysis,
  title={Analysis of optimal finite-element meshes in $R^1$},
  author={Babu{\v{s}}ka, Ivo and Rheinboldt, Werner C},
  journal={Mathematics of computation},
  volume={33},
  number={146},
  pages={435--463},
  year={1979}
}

@article{kelly1983posteriori,
  title={A posteriori error analysis and adaptive processes in the finite element method: Part I—error analysis},
  author={Kelly, Donald W and De SR Gago, JP and Zienkiewicz, Olgierd C and Babuska, Ivo},
  journal={International journal for numerical methods in engineering},
  volume={19},
  number={11},
  pages={1593--1619},
  year={1983},
  publisher={Wiley Online Library}
}

@book{bangerth2003adaptive,
  title={Adaptive finite element methods for differential equations},
  author={Bangerth, Wolfgang and Rannacher, Rolf},
  year={2003},
  publisher={Springer Science \& Business Media}
}

@article{ainsworth1997posteriori,
  title={A posteriori error estimation in finite element analysis},
  author={Ainsworth, Mark and Oden, J Tinsley},
  journal={Computer methods in applied mechanics and engineering},
  volume={142},
  number={1-2},
  pages={1--88},
  year={1997},
  publisher={Elsevier}
}

@article{becker2001optimal,
  title={An optimal control approach to a posteriori error estimation in finite element methods},
  author={Becker, Roland and Rannacher, Rolf},
  journal={Acta numerica},
  volume={10},
  pages={1--102},
  year={2001},
  publisher={Cambridge University Press}
}

@article{dorfler1996,
  title   = {A convergent adaptive algorithm for {Poisson’s} equation},
  author  = { Willy D{\"o}rfler},
  journal = {SIAM Journal on Numerical Analysis},
  volume  = {33},
  number  = {3},
  pages   = {1106--1124},
  year    = {1996},
  publisher = {SIAM}
}

@TechReport{      petsc-user-ref,
  author        = {Satish Balay and Shrirang Abhyankar and Mark~F. Adams and Steven Benson and Jed
                  Brown and Peter Brune and Kris Buschelman and Emil Constantinescu and Lisandro
                  Dalcin and Alp Dener and Victor Eijkhout and Jacob Faibussowitsch and William~D.
                  Gropp and V\'{a}clav Hapla and Tobin Isaac and Pierre Jolivet and Dmitry Karpeev
                  and Dinesh Kaushik and Matthew~G. Knepley and Fande Kong and Scott Kruger and
                  Dave~A. May and Lois Curfman McInnes and Richard Tran Mills and Lawrence Mitchell
                  and Todd Munson and Jose~E. Roman and Karl Rupp and Patrick Sanan and Jason Sarich
                  and Barry~F. Smith and Hansol Suh and Stefano Zampini and Hong Zhang and Hong Zhang
                  and Junchao Zhang},
  title         = {{PETSc/TAO} Users Manual},
  institution   = {Argonne National Laboratory},
  number        = {ANL-21/39 - Revision 3.22},
  doi           = {10.2172/2205494},
  year          = {2024}
}

@Misc{            petsc-web-page,
  author        = {Satish Balay and Shrirang Abhyankar and Mark~F. Adams and Steven Benson and Jed
                  Brown and Peter Brune and Kris Buschelman and Emil~M. Constantinescu and Lisandro
                  Dalcin and Alp Dener and Victor Eijkhout and Jacob Faibussowitsch and William~D.
                  Gropp and V\'{a}clav Hapla and Tobin Isaac and Pierre Jolivet and Dmitry Karpeev
                  and Dinesh Kaushik and Matthew~G. Knepley and Fande Kong and Scott Kruger and
                  Dave~A. May and Lois Curfman McInnes and Richard Tran Mills and Lawrence Mitchell
                  and Todd Munson and Jose~E. Roman and Karl Rupp and Patrick Sanan and Jason Sarich
                  and Barry~F. Smith and Stefano Zampini and Hong Zhang and Hong Zhang and Junchao
                  Zhang},
  title         = {{PETS}c {W}eb page},
  url           = {https://petsc.org/},
  howpublished  = {\url{https://petsc.org/}},
  year          = {2024}
}
\end{document}